\documentclass[leqno,12pt]{amsart}
\usepackage{amssymb} 
\theoremstyle{plain} 
\newtheorem{theorem}{\indent\sc Theorem}[section] 

\newtheorem{proposition}[theorem]{\indent\sc Proposition}

\theoremstyle{definition} 

\begin{document}

\title[Ordinary differential system]{Ordinary differential system in dinension six with affine Weyl group symmetry of type $D_4^{(2)}$ \\}

\author{Yusuke Sasano }

\renewcommand{\thefootnote}{\fnsymbol{footnote}}
\footnote[0]{2000\textit{ Mathematics Subjet Classification}.
34M55; 34M45; 58F05; 32S65.}

\maketitle

\begin{abstract}
We find a three-parameter family of ordinary differential systems in dimension six with affine Weyl group symmetry of type $D_4^{(2)}$.  This is the second example which gave higher order Painlev\'e type systems of type $D_{4}^{(2)}$. We show that we give its symmetry and holomorphy conditions. These symmetries, holomorphy conditions and invariant divisors are new.
\end{abstract}

\section{Introduction}
In \cite{Sasa10}, we study a two-parameter (resp. four-parameter) family of ordinary differential systems with affine Weyl group symmetry of type $D_3^{(2)}$ (resp. $D_5^{(2)}$). They are considered to be higher order versions of $P_{II}$. These systems are equevalent to the polynomial Hamiltonian systems, and can be considered to be 2-coupled (resp. 4-coupled) Painlev\'e II systems in dimension four (resp. eight).

We will complete the study of the above problem in a series of papers, for which this paper is the second, resulting in a series of equations for the remaining affine root systems of types $D_l^{(2)} \ (l=4,6,7,\ldots)$. This paper is the stage in this project where we find a 3-parameter family of ordinary differential systems in dimension six with affine Weyl group symmetry of type $D_4^{(2)}$ given by
\begin{align}\label{1}
\begin{split}
\frac{dx}{dt}&=\frac{\partial H}{\partial y}, \quad \frac{dy}{dt}=-\frac{\partial H}{\partial x}, \quad \frac{dz}{dt}=\frac{\partial H}{\partial w}, \quad \frac{dw}{dt}=-\frac{\partial H}{\partial z}, \quad \frac{dq}{dt}=\frac{\partial H}{\partial p}, \quad \frac{dp}{dt}=-\frac{\partial H}{\partial q}
\end{split}
\end{align}
with the polynomial Hamiltonian
\begin{align}\label{2}
\begin{split}
H =&\frac{1}{4t}y^3+\frac{3}{2}y^2+\frac{3\alpha_3-1}{t}xy+\frac{3}{4t}z^2w^2+\frac{3}{2}z^2w+\frac{3\alpha_1+3\alpha_2-2}{2t}zw+\frac{3}{2}\alpha_1 z\\
&-\frac{4}{t}p^3-6p^2-\frac{3\alpha_1+3\alpha_2+3\alpha_3-2}{t}qp-6tp+\frac{3}{4t}\alpha_1(8xp+2zp+yz)+\frac{6}{t}\alpha_2 xp\\
&+\frac{3}{2t}\alpha_3 (4xp-yq)+\frac{3}{4t}(8xyqp-4zwqp+8xzwp-8x^2yp+2z^2wp+yz^2w\\
&-2yq^2p+8w^2p-4yp^2+4yw^2+4y^2w+8ywp-8xp+8tyw).
\end{split}
\end{align}
Here $x,y,z,w,q$ and $p$ denote unknown complex variables, and $\alpha_0,\alpha_1,\alpha_2,\alpha_3$ are complex parameters satisfying the relation:
\begin{equation}\label{3}
\alpha_0+\alpha_1+\alpha_2+\alpha_3=1.
\end{equation}
In section 2, each principal part of this Hamiltonian can be transformed into the one with its first integrals by birational and symplectic transformations. However, the Hamiltonian $H$ is not the first integral.

We remark that for this system we tried to seek its first integrals of polynomial type with respect to $x,y,z,w,q,p$. However, we can not find.

This is the second example which gave higher order Painlev\'e type systems of type $D_{4}^{(2)}$.

We also remark that 2-coupled Painlev\'e III system in dimension four given in the paper \cite{Sasa2} admits the affine Weyl group symmetry of type $D_4^{(2)}$ as the group of its B{\"a}cklund transformations, whose generators $w_1,w_2$ are determined by the invariant divisors. However, the transformations $w_3,w_4$ do not satisfy so (see Theorem 4.1 in \cite{Sasa2}).

On the other hand, the system \eqref{1} admits the affine Weyl group symmetry of type $D_4^{(2)}$ as the group of its B{\"a}cklund transformations, whose generators $s_0,\ldots,s_3$ are determined by the invariant divisors \eqref{invariant}.

\section{Principal parts of the Hamiltonian}
In this section, we study three Hamiltonians $K_1,K_2$ and $K_3$ in the Hamiltonian $H$.

At first, we study the Hamiltonian system
\begin{align}
\begin{split}
\frac{dx}{dt}=&\frac{\partial K_1}{\partial y}=\frac{3y(y+4t)-4(\alpha_0+\alpha_1+\alpha_2-2\alpha_3)x}{4t},\\
\frac{dy}{dt}=&-\frac{\partial K_1}{\partial x}=\frac{(\alpha_0+\alpha_1+\alpha_2-2\alpha_3)y}{t}
\end{split}
\end{align}
with the polynomial Hamiltonian
\begin{align}\label{4}
\begin{split}
K_1 =&\frac{1}{4t}y^3+\frac{3}{2}y^2+\frac{3\alpha_3-1}{t}xy,
\end{split}
\end{align}
where setting $z=w=q=p=0$ in the Hamiltonian $H$, we obtain $K_1$.

This equation can be explicitly solved by
\begin{align}
\begin{split}
x(t)=&\frac{C_1 t^{-1+3(\alpha_0+\alpha_1+\alpha_2)}}{\alpha_0+\alpha_1+\alpha_2-\alpha_3}+\frac{C_1^2 t^{-4+6(\alpha_0+\alpha_1+\alpha_2)}}{4(\alpha_0+\alpha_1+\alpha_2-2\alpha_3)}+C_2 t^{2-3(\alpha_0+\alpha_1+\alpha_2)},\\
y(t)=&C_1 t^{(\alpha_0+\alpha_1+\alpha_2-2\alpha_3)} \quad (C_1,C_2:integral \ constants).
\end{split}
\end{align}

Next, we study the Hamiltonian system
\begin{align}
\begin{split}
\frac{dz}{dt}&=\frac{\partial K_2}{\partial w}, \quad \frac{dw}{dt}=-\frac{\partial K_2}{\partial z}
\end{split}
\end{align}
with the polynomial Hamiltonian
\begin{align}\label{5}
\begin{split}
K_2 =&\frac{3}{4t}z^2w^2+\frac{3}{2}z^2w+\frac{3\alpha_1+3\alpha_2-2}{2t}zw+\frac{3}{2}\alpha_1 z,
\end{split}
\end{align}
where setting $x=y=q=p=0$ in the Hamiltonian $H$, we obtain $K_2$.

{\bf Step 1:} We make the change of variables:
\begin{equation}
z_1=tz, \quad w_1=\frac{w}{t}.
\end{equation}
We remark that this transformation is symplectic.

It is easy to see that the system with the polynomial Hamiltonian
\begin{equation}
\tilde{K}_2=\frac{3z_1(z_1w_1^2+2z_1w_1+2(\alpha_1+\alpha_2)w_1+2\alpha_1}{4t}
\end{equation}
has its first integral I:
\begin{equation}
I:=4t \tilde{K}_2.
\end{equation}

Finally, we study the Hamiltonian system
\begin{align}
\begin{split}
\frac{dq}{dt}=&\frac{\partial K_3}{\partial p}=-\frac{12p(p+t)-(2\alpha_0-\alpha_1-\alpha_2-\alpha_3)q+6t^2}{t},\\
\frac{dp}{dt}=&-\frac{\partial K_3}{\partial q}=-\frac{(2\alpha_0-\alpha_1-\alpha_2-\alpha_3)p}{t}
\end{split}
\end{align}
with the polynomial Hamiltonian
\begin{align}\label{6}
\begin{split}
K_3 =&-\frac{4}{t}p^3-6p^2-\frac{3\alpha_1+3\alpha_2+3\alpha_3-2}{t}qp-6tp,
\end{split}
\end{align}
where setting $x=y=z=w=0$ in the Hamiltonian $H$, we obtain $K_3$.

This equation can be explicitly solved by
\begin{align}
\begin{split}
q(t)=&-2t^2 \left\{\frac{1}{\alpha_1+\alpha_2+\alpha_3}+2C_1 t^{-6\alpha_0} \left(\frac{t^{3\alpha_0}}{-\alpha_0+\alpha_1+\alpha_2+\alpha_3}+\frac{C_1}{-2\alpha_0+\alpha_1+\alpha_2+\alpha_3} \right) \right\}\\
&+C_2t^{-1+3\alpha_0},\\
p(t)=&C_1 t^{-(2\alpha_0-\alpha_1-\alpha_2-\alpha_3)} \quad (C_1,C_2:integral \ constants).
\end{split}
\end{align}

\section{Symmetry and holomorphy conditions}
In this section, we study the symmetry and holomorphy conditions of the system \eqref{1}. These properties are new.

\begin{theorem}\label{th:1}
The system \eqref{1} admits the affine Weyl group symmetry of type $D_4^{(2)}$ as the group of its B{\"a}cklund transformations, whose generators $s_0,s_1,\ldots,s_3$ defined as follows$:$ with {\it the notation} $(*):=(x,y,z,w,q,p,t;\alpha_0,\alpha_1,\alpha_2,\alpha_3)$\rm{: \rm}
\begin{align}
\begin{split}
s_0:(*) \rightarrow &\left(x+\frac{\alpha_0}{y+z^2/4},y,z,w-\frac{\alpha_0 z}{2(y+z^2/4)},q,p,t;-\alpha_0,\alpha_1+2\alpha_0,\alpha_2,\alpha_3 \right),\\
s_1:(*) \rightarrow &\left(x,y,z+\frac{\alpha_1}{w},w,q,p,t;\alpha_0+\alpha_1,-\alpha_1,\alpha_2+\alpha_1,\alpha_3 \right),\\
s_2:(*) \rightarrow &(x+\frac{\alpha_2/2}{f_2},y-\frac{\alpha_2 z/2}{f_2}-\frac{\alpha_2^2/4}{f_2^2},z+\frac{\alpha_2}{f_2},w+\frac{\alpha_2(q-2x)/4}{f_2},\\
&q+\frac{\alpha_2}{f_2},p+\frac{\alpha_2 z/4}{f_2}+\frac{\alpha_2^2/8}{f_2^2},t;\alpha_0,\alpha_1+\alpha_2,-\alpha_2,\alpha_3+\alpha_2),\\
s_3:(*) \rightarrow &\left(x,y,z,w,q+\frac{\alpha_3}{p},p,t;\alpha_0,\alpha_1,\alpha_2+2\alpha_3,-\alpha_3 \right),
\end{split}
\end{align}
where $f_2:=w+p+\frac{y}{2}+\frac{xz}{2}-\frac{zq}{4}+t$.
\end{theorem}
We note that the B{\"a}cklund transformations of this system satisfy
\begin{equation}
s_i(g)=g+\frac{\alpha_i}{f_i}\{f_i,g\}+\frac{1}{2!} \left(\frac{\alpha_i}{f_i} \right)^2 \{f_i,\{f_i,g\} \}+\cdots \quad (g \in {\Bbb C}(t)[x,y,z,w,q,p]),
\end{equation}
where poisson bracket $\{,\}$ satisfies the relations:
$$
\{y,x\}=\{w,z\}=\{p,q\}=1, \quad the \ others \ are \ 0.
$$
Since these B{\"a}cklund transformations have Lie theoretic origin, similarity reduction of a Drinfeld-Sokolov hierarchy admits such a B{\"a}cklund symmetry.

\begin{proposition}
This system has the following invariant divisors\rm{:\rm}
\begin{center}\label{invariant}
\begin{tabular}{|c|c|c|} \hline
parameter's relation & $f_i$ \\ \hline
$\alpha_0=0$ & $f_0:=y+\frac{z^2}{4}$  \\ \hline
$\alpha_1=0$ & $f_1:=w$  \\ \hline
$\alpha_2=0$ & $f_2:=w+p+\frac{y}{2}+\frac{xz}{2}-\frac{zq}{4}+t$  \\ \hline
$\alpha_3=0$ & $f_3:=p$  \\ \hline
\end{tabular}
\end{center}
\end{proposition}
We note that when $\alpha_1=0$, we see that the system \eqref{1} admits a particular solution $w=0$, and when $\alpha_2=0$, after we make the birational and symplectic transformations:
\begin{equation}
x_2=x-\frac{z}{2}, \ y_2=y+\frac{z^2}{4}, \ z_2=z, \ w_2=w+\frac{y}{2}+p+t+\frac{xz}{2}-\frac{zq}{4}, \ q_2=q-z, \quad p_2=p-\frac{z^2}{8}.
\end{equation}
we see that the system \eqref{1} admits a particular solution $w_2=0$.

\begin{proposition}
Let us define the following translation operators{\rm : \rm}
\begin{align}
\begin{split}
&T_1:=s_1 s_2 s_3 s_2 s_1 s_0, \quad T_2:=s_1 T_1 s_1, \quad T_3:=s_2 T_2 s_2.
\end{split}
\end{align}
These translation operators act on parameters $\alpha_i$ as follows$:$
\begin{align}
\begin{split}
T_1(\alpha_0,\alpha_1,\alpha_2,\alpha_3)=&(\alpha_0,\alpha_1,\alpha_2,\alpha_3)+(-2,2,0,0),\\
T_2(\alpha_0,\alpha_1,\alpha_2,\alpha_3)=&(\alpha_0,\alpha_1,\alpha_2,\alpha_3)+(0,-2,2,0),\\
T_3(\alpha_0,\alpha_1,\alpha_2,\alpha_3)=&(\alpha_0,\alpha_1,\alpha_2,\alpha_3)+(0,0,-2,2).
\end{split}
\end{align}
\end{proposition}

\begin{theorem}\label{pro:2}
Let us consider a polynomial Hamiltonian system with Hamiltonian $K \in {\Bbb C}(t)[x,y,z,w,q,p]$. We assume that

$(A1)$ $deg(K)=4$ with respect to $x,y,z,w,q,p$.

$(A2)$ This system becomes again a polynomial Hamiltonian system in each coordinate system $r_i \ (i=0,1,2,3)${\rm : \rm}
\begin{align}
\begin{split}
r_0:&x_0=\frac{1}{x}, \ y_0=-\left( \left(y+\frac{z^2}{4} \right)x+\alpha_0 \right)x, \ z_0=z, \quad w_0=w+\frac{xz}{2}, \ q_0=q, \quad p_0=p,\\
r_1:&x_1=x, \ y_1=y, \ z_1=\frac{1}{z}, \ w_1=-(wz+\alpha_1)z, \ q_1=q, \quad p_1=p,\\
r_2:&x_2=x-\frac{z}{2}, \ y_2=y+\frac{z^2}{4}, \ z_2=\frac{1}{z}, \ w_2=-\left( \left(w+\frac{y}{2}+p+t+\frac{xz}{2}-\frac{zq}{4} \right)z+\alpha_2 \right)z,\\
&q_2=q-z, \quad p_2=p-\frac{z^2}{8}, \\
r_3:&x_3=x, \ y_3=y, \ z_3=z, \ w_3=w, \ q_3=\frac{1}{q}, \quad p_3=-(pq+\alpha_3)q.
\end{split}
\end{align}
Then such a system coincides with the system
\begin{align}\label{8}
\begin{split}
\frac{dx}{dt}&=\frac{\partial K}{\partial y}, \quad \frac{dy}{dt}=-\frac{\partial K}{\partial x}, \quad \frac{dz}{dt}=\frac{\partial K}{\partial w}, \quad \frac{dw}{dt}=-\frac{\partial K}{\partial z}, \quad \frac{dq}{dt}=\frac{\partial K}{\partial p}, \quad \frac{dp}{dt}=-\frac{\partial K}{\partial q}
\end{split}
\end{align}
with the polynomial Hamiltonian
\begin{align}\label{9}
\begin{split}
K =&H+a_1(y+2p)+a_2(y+2p)^2+a_3(y+2p)^3+a_4(y+2p)^4 \quad (a_i \in {\Bbb C}(t)).
\end{split}
\end{align}
\end{theorem}
We note that the condition $(A2)$ should be read that
\begin{align*}
&r_j(K) \quad (j=0,1,3), \quad r_2(K-z)
\end{align*}
are polynomials with respect to $x_i,y_i,z_i,w_i,q_i,p_i$.

We remark that $y+2p$ is not the first integral of the system \eqref{8} with the polynomial Hamiltonian \eqref{9}.

\end{document}